\documentclass[twoside]{article}
\usepackage{amsmath,amsthm,amssymb}
\usepackage{url}
\usepackage[spanish,english]{babel}
\usepackage{fancyhdr}
\fancypagestyle{plain}{\fancyhf{}%
\fancyhead[L]{\footnotesize{\pageref{begin-art}--\pageref{end-art}}}%
}

\fancypagestyle{headings}{\fancyhf{}%
\fancyhead[LE,RO]{\thepage}%
\fancyfoot[C]{\footnotesize{\pageref{begin-art}--\pageref{end-art}}}}

\def \R {\mathbb R}
\def \C {\mathbb C}
\def \G {\mathcal G}
\def \F {\mathcal F}
\def \Fo {\mathcal F_o}

\def \y {{\mathrm y}}
\def \and {{\mathrm and}}
\def \R {\mathbb R}
\def \C {\mathbb C}
\def \G {\mathcal G}
\def \F {\mathcal F}
\def \Fo {\mathcal F_o}

\def \d {\partial}
\def \yu {{\underbar y}}
\def \Fu {{\bf F}_1}
\def \Fd {{\bf F}_2}
\def \y {{\mathrm {and}}}
\def \yb {{\bar y}}
\def \yu {{\underbar y}}
\def \xk {x_{\hat k}}
\def \yk {y_{\hat k}}
\def \yyk {y'_{\hat k}}
\def \yuk {\yu_{\hat k}}
\def \u {{\underbar 1}}

\begin{document}
\label{begin-art}
\pagestyle{headings}
\thispagestyle{plain} 
\footnote{\hspace{-18.1pt} }
\selectlanguage{english}
\begin{center}
{\large\bfseries The Clifford-Fourier Transform $\Fo$ and Monogenic Extensions \par }

\vspace{5mm}
{\large Arnoldo Bezanilla L\'opez (\url{abeza@fcfm.buap.mx})}\\[1mm]
Benem\'erita Universidad Aut\'onoma de Puebla \\
72570. Puebla, Puebla, M\'exico \\

\

{\large Omar Le\'on S\'anchez (\url{oleonsan@math.uwaterloo.ca})}\\[1mm]
University of Waterloo\\
N2L 3G1. Waterloo, Ontario, Canada\\

\end{center}

\selectlanguage{english}

\begin{abstract} In the last decade several versions of the Fourier transform have been formulated in the framework of Clifford algebra. We present a (Clifford-Fourier) transform, constructed using the geometric properties of Clifford algebra. We show the corresponding results of operational calculus, and a connection between the Fourier transform and this new transform. We obtain a technique to construct monogenic extensions of a certain type of continuous functions, and versions of the Paley-Wiener theorems are formulated.

\

\end{abstract}

\selectlanguage{english}

\section{Introduction}

It is well known that the convolution and the Fourier transform are two commonly used tools in image processing of vector fields. Due to several important applications, these techniques have been extended to the analysis of multivector fields and geometric algebras \cite{BaySob}. In \cite{Bu}, B\"ulow defines a Clifford-Fourier transform by
\begin{displaymath}\label{ClifFo}
(\F_b f)(y)=\int_{\R^n}f(x)e^{-2\pi e_1x_1y_1}\cdots e^{-2\pi e_nx_ny_n}dx
\end{displaymath}
where the product in the kernel of the transform must be performed in the order determined by the indices. This makes the convolution theorem look rather complicated. This kernel had already been introduced as a theorical concept in Clifford analysis \cite{Som}, and Sommen studied it in connection with generalizations of the Laplace, Cauchy and Hilbert transform \cite{Som2}.
In \cite{Eb}, Ebling and Scheuermann study the convolution of their Clifford-Fourier transforms: for $f\colon\R^2 \to \G_2$
\begin{displaymath}
(\F_e f)(y)=\frac{1}{2\pi}\int_{\R^2}e^{- e_{12}x\cdot y}f(x)dx
\end{displaymath}
and for $f\colon\R^3\to \G_3$
\begin{displaymath}
(\F_e f)(y)=\frac{1}{(2\pi)^{\frac{3}{2}}} \int_{\R^3}e^{- e_{123}x\cdot y}f(x)dx.
\end{displaymath}
These transforms allowed them to prove the convolution and differentiation theorems as extensions of the usual results on scalar fields. However, there is no similar development for higher dimension spaces. For more on Clifford-Fourier transforms see \cite{Bra}, \cite{Bra2}, \cite{Bra3}.

We introduce a Clifford-Fourier transform, using the basic properties of the Clifford algebra. We use the identification $\G_n\subseteq\G_{2n}$ to build $n$ distinct bivectors, where each of these will behave as an imaginary unit. Hence, the Clifford-Fourier transform will not be considered as a generalization of the Fourier transform, but as a formulation of this in the context of Clifford algebra.

It is known that the Fourier transform can be used to construct analytic functions in $\C^n$, and that the Paley-Wiener theorems tell us when the restriction to $\R^n$ of an analytic function is the Fourier transform of a square integrable function \cite{St}. In \cite{KitTao}, Kou and Qian give a proof of the Paley-Wiener theorems for the Fourier transform in the Clifford algebra setting. Similarly, we use the Clifford-Fourier transform to obtain monogenic functions in $\R^{2n}$ and versions of the Paley-Wiener theorems.

This paper is mostly a presentation of the basic properties of the Clifford-Fourier transform, and it presents the advantages that this particular transform has in the construction of monogenic functions. In section 2, we present the basic definitions of the Clifford algebra and monogenic functions. In section 3, we define the Clifford-Fourier transform and develop the operational calculus. Finally, in section 4, we construct monogenic extensions using the Clifford-Fourier transform, and we prove versions of the Paley-Wiener theorems.

\

\section{Preliminaries}

\subsection{The Clifford Algebra}

Consider the Euclidean space $\R^n$, with the usual scalar product
\begin{displaymath}
x\cdot y=\sum_{i=1}^n x_iy_i
\end{displaymath}
The Clifford algebra $\G_n$ is defined as the associative algebra with unity over $\R$, that contains $\R$ and $\R^n$ as distinct subspaces; is generated by $\R^n$, but not by any proper subspace and
\begin{displaymath}
x^2=x\cdot x, \qquad \mathrm{for \ all} \ x\in \R^n.
\end{displaymath}

Using an orthonormal basis $\{e_1,\dots,e_n\}$ of $\R^n$, one can give a formal construction of the Clifford algebra. That is, each multivector $A\in \G_n$ is represented in the form
\begin{displaymath}
A=\sum_{I}\alpha_Ie_I, \quad \alpha_I \in \R,
\end{displaymath}
where $e_I=e_{i_1}\cdots e_{i_k}$ with $i_1<\cdots<i_k$ and $e_{\emptyset}=1$. Multiplication is determined by
\begin{displaymath}\label{regla}
e_i^2=1 \qquad \and \qquad e_ie_j=-e_je_i, \quad i\neq j.
\end{displaymath}
Any multivector can be expressed in the form
\begin{displaymath}
A=\sum_{k=0}^n<A>_k, 
\end{displaymath}
where
\begin{displaymath}
<A>_k=\sum_{|I|=k}\alpha_Ie_I,
\end{displaymath}
the last expression is known as the $k$-vector component of $A$.

The reversion is the anti-involution ${}^{\dag}\colon\G_n\to\G_n$ such that $x^{\dag}=x$ for any vector $x\in\R^n$. Also, the scalar product of multivectors is defined by
\begin{equation}\label{prodesc}
A\cdot B=<A\,B^{\dag}>_0
\end{equation}
In particular, one has that if $A=\sum_{I}\alpha_Ie_I$
\begin{equation}\label{norma}
|A|^2= A\cdot A=<A\,A^{\dag}>_0=\sum_I\alpha_I^2.
\end{equation}
This defines a norm in $\G_n$ with the property
\begin{displaymath}
|AB|\leq 2^{\frac{n}{2}}|A||B|
\end{displaymath}
and if $AA^{\dag}=<AA^{\dag}>_0$, then $|AB|=|A||B|$. Hence, $\G_n$ is a normed associative non-commutative algebra with unity.

For a set of multivectors $M \subseteq \G_n$, we will use the notation $\G(M)$ for the subalgebra of $\G_n$ generated by $M$. In particular, $\G(e_1,\dots,e_n)=\G_n$.

\

\subsection{Monogenic Functions}

The vector derivative in $\R^n$ is defined by
\begin{displaymath}
\partial =\sum_{i=1}^n e_i\frac{\partial }{\partial x_i}
\end{displaymath}
If $f\colon\R^n\to\G_{n}$ is of class $C^1$, we say that $f$ is monogenic (or left-monogenic) in a domain $M\subseteq \R^n$ if
\begin{displaymath}
\partial f(x)=0, \quad  \textrm{for all} \quad x\in  M,
\end{displaymath}
and right-monogenic if $f(x)\partial=0$.

\

\noindent {\bf Theorem 2.1.} Let $M$ be a piecewise differentiable bounded domain in $\R^n$ with boundary $\beta(M)$. If $f,g\colon \R^n\to \G_n$ are of class $C^1$, then
\begin{displaymath}
\int_M \d f(x) dx=\int_{\beta(M)}n(x)f(x) dx
\end{displaymath}
and 
\begin{displaymath}
\int_M g(x)\d dx=\int_{\beta(M)}g(x)n(x) dx
\end{displaymath}
where $n(x)$ is the outer normal to $\beta(M)$ in $x$. 

\

This theorem is known as the fundamental theorem of Clifford's geometric calculus \cite{He}. As a consequence one has Clifford-Cauchy's theorem:

\

\noindent {\bf Theorem 2.2.} Let $M$, $f$ and $g$ as in theorem 2.1. If $f$ is monogenic and $g$ is right-monogenic in $M$, then
\begin{displaymath}
\int_{\beta(M)}n(x)f(x) dx=0
\end{displaymath}
and
\begin{displaymath}
\int_{\beta(M)}g(x)n(x) dx=0
\end{displaymath}

\

Another consequence is Clifford-Cauchy's integral formula:

\

\noindent {\bf Theorem 2.3.} Let $M$, $f$ and $g$ as in theorem 2.1. If $f$ is monogenic and $g$ is right-monogenic in $M$ and $x\in M$ is a point in the interior, then
\begin{displaymath}
f(x)=\frac{1}{\Omega_n}\int_{\beta(M)}\frac{y-x}{|y-x|^n}n(y)f(y) dy
\end{displaymath}
and
\begin{displaymath}
g(x)=\frac{1}{\Omega_n}\int_{\beta(M)}g(y)n(y)\frac{y-x}{|y-x|^n} dy
\end{displaymath}
where $\Omega_n$ is the area of the unit $(n-1)$-sphere in $\R^n$.

\

Because of these theorems (and many others), monogenic functions are considered a generalization of complex analytic functions \cite{Som}, \cite{He}, \cite{Ab}.

\

\section{The Clifford-Fourier transform $\Fo$}

\subsection{Definition}

All integrals will be considered with respect to the normalized measure:
\begin{displaymath}
dm(x)=(2 \pi)^{-\frac{n}{2}}dx,
\end{displaymath}
where $dx$ is the Lebesgue measure in $\R^n$. The norms in $L_1$ and $L_2$ will be taken with respect to this measure, i.e.
\begin{displaymath}
\|f\|_p=\left(\int_{\R^n}|f(x)|^{p}dm(x)\right)^{\frac{1}{p}},\quad p=1,2,
\end{displaymath}
with $f$ taking its values in the Clifford algebra and $|f(x)|$ being the norm defined in (\ref{norma}).

Let 
\begin{displaymath}
\{e_1,\dots,e_n,e'_1,\dots, e'_n\}\subseteq \R^{2n}
\end{displaymath}
be an orthonormal basis of $\R^{2n}$, and $\G_{2n}$ the Clifford algebra generated by these vectors. Identifying $\R^n$ with the subspace generated by $\{e_1,\dots,e_n\}$, we consider $\R^n\subseteq \R^{2n}$ and $\G_n\subseteq \G_{2n}$.

Define the bivectors
\begin{displaymath}
B_i=e_ie'_i, \qquad i=1,2,\dots,n,
\end{displaymath}
one has that
\begin{displaymath}
B_i^2=-1 \qquad \textrm{and} \qquad B_iB_j=B_jB_i,\qquad i,j=1,\dots,n.
\end{displaymath}
Note that $\G(B_1,\dots,B_n)$ is a commutative subalgebra of $\G_{2n}$.

Now, consider the bilinear symmetric function $I\colon\R^n\times\R^n\to \G_{2n}$
\begin{displaymath}
I(x,y)=\sum_{i=1}^n B_iP_i(x)P_i(y)
\end{displaymath}
where $P_i$ denotes the projection of $\R^{2n}$ to the subspace generated by $\{e_i,e'_i\}$. Note that the function $I$ can be written as
\begin{displaymath}
I(x,y)=B_1x_1y_1+\cdots+B_nx_ny_n.
\end{displaymath}
The reason for using the projection $P_i$ in the definition of $I(x,y)$ will be clarified in section 4 when we generalize the Clifford-Fourier transform for the construction of monogenic functions.

To define the Clifford-Fourier transform, we use the following kernel
\begin{equation}\label{nucleo}
e^{-I(x,y)}=e^{-\sum B_ix_iy_i}=\prod_{i=1}^n e^{-B_ix_iy_i}
\end{equation}
and note that
\begin{displaymath}
|e^{-I(x,y)}|=1.
\end{displaymath}

\

\noindent {\bf Definition 3.1.} Let $f\colon\R^n\to\G_{2n}$, assuming the integrals exist, we define the left Clifford-Fourier transform as
\begin{equation}\label{defClifTrans1}
(\Fo f)(y)=\int_{\R^n}e^{-I(x,y)}f(x)dm(x)
\end{equation}
and the right Clifford-Fourier transform as
\begin{equation}\label{defizq}
(f \Fo)(y)=\int_{\R^n}f(x)e^{-I(x,y)}dm(x).
\end{equation}

Note that if $Im(f)\subseteq \G(I_1,\dots,I_n)$, then
\begin{displaymath}
\Fo f=f \Fo.
\end{displaymath}

\

\subsection{Connection between $\F$ and $\Fo$}

For the purposes of this section we extend the range of the functions from $\G_{2n}$ to $\C\otimes\G_{2n}$, the complexification of $\G_{2n}$ with the tensor product of algebras. For $f\colon \R^n\to \C\otimes\G_{2n}$, the Fourier transform can be defined by
\begin{displaymath}
(\F f)(y)=\int_{\R^n}e^{-ix\cdot y}f(x)dm(x).
\end{displaymath}

We will use the following multivectors
\begin{displaymath}
M^+_k=\frac{1}{2}(1+iB_k),\qquad M^-_k=\frac{1}{2}(1-iB_k),\qquad k=1,\dots,n.
\end{displaymath}
They safisfy
\begin{itemize}
\item $M^+_k+M^-_k=1$, \quad $M^+_kM^-_k=0$, \quad $(M^{\pm}_k)^2=M^{\pm}_k$.
\item $M^+_k B_k=M^+_k (-i)$, \quad $M^-_kB_k=M^-_ki$.
\end{itemize}

\

Let $s\in{\cal J}= \{-1,1\}^n$, then $s=(s_1,\dots,s_n)$ where $s_k=1$ or $s_k=-1$. For each $s$ the reflection in $\R^n$ is defined by:
\begin{displaymath}
R_s(y)=-(s_1y_1e_1+\cdots+s_ny_ne_n).
\end{displaymath}

Defining $M^{s_k}_k=M^{\pm}_k$ according to $s_k=\pm 1$, and $M^s=M^{s_1}_1 M^{s_2}_2\cdots M^{s_n}_n$. We have
\begin{displaymath}
\sum_{s\in {\cal J}}M^s=1, \qquad M^sM^{s'}=0 \quad (s\neq s'), \qquad (M^s)^2=M^s.
\end{displaymath}
Hence
\begin{displaymath}
e^{-ix\cdot y}=\prod_{k=1}^n(M^+_k + M^-_k)e^{-ix\cdot y}=\sum_{s\in {\cal J}}M^s e^{-I(x,R_s(y))}
\end{displaymath}
Therefore, we have the following relation
\begin{displaymath}
(\F f)(y)=\sum_{s\in {\cal J}}M^s (\Fo f)(R_s(y)).
\end{displaymath}

From this we get
\begin{eqnarray*}
(\F(\F f))(x)
&=& \sum_{r\in {\cal J}_n}\sum_{s\in {\cal J}_n}M^rM^s(\Fo(\Fo f(R_s(y))))(R_r(x))\\
&=& \sum_{s\in {\cal J}_n}M^s(\Fo(\Fo f))(R_s^2(x))\\
&=& \sum_{s\in {\cal J}_n}M^s(\Fo(\Fo f))(x) \\
&=& (\Fo(\Fo f))(x) \\
\end{eqnarray*}
that is, $\F^2 f=\Fo^2 f$.

\

On the other hand, one can get a similar equation for $\Fo f$ in terms of $\F f$:
\begin{displaymath}
(\Fo f)(y)=\sum_{s\in {\cal J}}M^s (\F f)(R_s(y)).
\end{displaymath}

In an analogous way, we obtain the following expressions
\begin{displaymath}
(\F f )(y)=\sum_{s\in {\cal J}}(f\Fo )(R_s(y))M^s 
\end{displaymath}
and
\begin{displaymath}
(f\Fo )(y)=\sum_{s\in {\cal J}}(\F f)(R_s(y))M^s .
\end{displaymath}

\

\subsection{Operational Calculus of $\Fo$}

For each integrable function $f\colon\R^n\to\G_{2n}$, the transforms $\Fo f$ and $f\Fo$ are well defined, belong to $L_{\infty}$ and
\begin{displaymath}
\|\Fo f\|_{\infty}\leq 2^n\|f\|_1  \qquad \and \qquad \|f\Fo \|_{\infty}\leq 2^n\|f\|_1
\end{displaymath}
hence, $f\to \Fo f$ and $f\to f\Fo$ are bounded maps from $L_1\to L_{\infty}$.

\

\noindent {\bf Lemma 3.1.} Let $f$, $g\in L_1$.
\begin{enumerate}
\item If $A$ and $B\in \G_{2n}$, 
\begin{displaymath}
\Fo(fA+gB)=(\Fo f)A+(\Fo g)B \quad \text{and} \quad (Af+Bg)\Fo=A(f\Fo)+B(g\Fo).
\end{displaymath}
\item If $u\in \R^n$ and $\tau_u f(x):=f(x-u)$,
\begin{displaymath}
(\Fo(\tau_{u}f))(y)=e^{-I(u,y)}(\Fo f)(y) \quad \text{and} \quad ((\tau_{u}f)\Fo)(y)=(f\Fo)(y)e^{-I(u,y)}.
\end{displaymath}
\item If $h_1(x)=e^{I(x,u)}f(x)$ and $h_2(x)=f(x)e^{I(x,u)}$,
\begin{displaymath}
\Fo h_1=\tau_{u} (\Fo f) \quad \text{and} \quad h_2\Fo =\tau_{u} (f\Fo).
\end{displaymath}
\item If $\lambda >0$ and $h(x)=f(x/ \lambda)$,
\begin{displaymath}
(\Fo h)(y)=\lambda^n (\Fo f)(\lambda y) \quad \text{and} \quad (h\Fo )(y)=\lambda^n (f\Fo )(\lambda y).
\end{displaymath}
\item If $h(x)=f^{\dag}(-x)$,
\begin{displaymath}
\Fo h=(f \Fo)^{\dag} \quad \text{and} \quad  h\Fo=(\Fo f )^{\dag}.
\end{displaymath}
\end{enumerate}

\

Recall that the convolution of the integrable functions $f$ and $g$ is defined almost everywhere by
\begin{displaymath}
f*g \ (y)=\int_{\R^n}f(y-x)g(x)dm(x).
\end{displaymath}

\

\noindent {\bf Convolution Theorem.} Let $f$ and $g\in L_1$. If $Im(f)\subseteq \G(I_1,\dots,I_n)$, then
\begin{displaymath}
\Fo(f*g)=(\Fo f)(\Fo g).
\end{displaymath}
If $Im(g)\subseteq \G(I_1,\dots,I_n)$, then
\begin{displaymath}
(f*g)\Fo=(f \Fo )(g\Fo ).
\end{displaymath}
{\bf Proof.}
If $Im(f)\subseteq \G_I$, $f$ commutes with the kernel of the transform and applying Fubini's theorem
\begin{displaymath}
\Fo (f* g)(y)=\int_{\R^n}\int_{\R^n}e^{-I(x,y)}f(x-t)g(t)dm(t)dm(x) 
\end{displaymath}
\begin{displaymath}
=\int_{\R^n}\int_{\R^n}e^{-I(x-t,y)}f(x-t)e^{-I(t,y)}g(t)dm(t)dm(x)
\end{displaymath}
\begin{displaymath}
= \int_{\R^n}e^{-I(x,y)}f(x)dm(x)\int_{\R^n}e^{-I(t,y)}g(t)dm(t)=(\Fo f)(y)(\Fo g)(y).  
\end{displaymath}

The other case is analogous.

\hspace*{11cm} $\Box$

\

In analogy with \cite{Eb}, these formulations of the convolution theorem for the Clifford-Fourier transform suggest certain applications in computer graphics.

\

\noindent {\bf Multiplication Formula.} If $f$ and $g\in L_1$, then
\begin{displaymath}
\int_{\R^n}(f\Fo)(x)g(x)dm(x)=\int_{\R^n}f(x)(\Fo g)(x)dm(x)
\end{displaymath}

\

The image of the Clifford-Fourier transform of an integrable function is a continuous function that vanishes at infinity, i.e. if $C_0$ denotes the space of continuous functions defined in $\R^n$ with values in $\G_{2n}$ that vanish at infinity, one has

\

\noindent {\bf Riemann-Lebesgue theorem.} If $f\in L_1$, then $\Fo f$ and $f \Fo \in C_0$.

\

\noindent {\bf Proof.} The continuity follows from the dominated convergence theorem. Let $y=y_1e_1+\cdots +y_ne_n$ and suppose that some $y_k\neq 0$. Define $y_k^*=(\pi/y_k)e_k$, then
\begin{displaymath}
|(\Fo f)(y)|\leq {2^{n-1}}\|f-(\tau_{y_k^*}f)\|_1
\end{displaymath}
On the other hand, if $\epsilon >0$ one can find $\delta >0$ such that if $|y|>\delta$ then
\begin{displaymath}
\|f-(\tau_{y_k^*}f)\|_1< \frac{\epsilon}{2^{n-1}} \; 
\end{displaymath}
for some $k$. Hence, if $|y|>\delta$, then $|(\Fo f)(y)|<\epsilon$. The proof for $f\Fo$ is analogous.

\hspace*{11cm} $\Box$

\ 

\noindent {\bf Differentiation theorem.} Let $f\in L_1$.
\begin{itemize}
\item Let $h_1(x)=-B_k x_kf(x)$ and $h_2(x)=- f(x)B_kx_k$. 

If $h_1\in L_1$, 
\begin{displaymath}
\frac{\partial(\Fo f)}{\partial y_k}=\Fo h_1.
\end{displaymath}
If $h_2\in L_1$,
\begin{displaymath}
\frac{\partial(f\Fo )}{\partial y_k}=h_2\Fo .
\end{displaymath}
\item If $\frac{\partial f}{\partial x_k}(x)$ exists for almost all $x$ and is integrable, then 
\begin{displaymath}
\left(\Fo \frac{\partial f}{\partial x_k}\right)(y)=B_k y_k(\Fo f)(y)
\end{displaymath}
and
\begin{displaymath}
\left(\frac{\partial f}{\partial x_k}\Fo \right)(y)=(f\Fo )(y)B_k y_k.
\end{displaymath}
\end{itemize}

\

Hence, in analogy with the Fourier transform, the Clifford-Fourier transform can be used as a tool to study several aspects of applied problems, such as multivector differential equations.

\

\noindent {\bf Inversion Theorem.} Let $f\in L_1$.
\begin{itemize}
\item If $\Fo f \in L_1$,
\begin{displaymath}
f(x)=(\Fo^2 f)(-x) \quad \text{almost everywhere}.
\end{displaymath}
\item If $f \Fo \in L_1$,
\begin{displaymath}
f(x)=( f \Fo^2)(-x) \quad \text{almost everywhere}.
\end{displaymath}
\end{itemize}
Moreover, one has that $\Fo^2 f=f \Fo^2$.

\

On the other hand, for $f\in L_2$ its Clifford-Fourier transform is not always well defined. However, if $f\in L_1\cap L_2$ one has that $\Fo f$ and $f \Fo \in L_2$. Moreover, we have the following proposition:

\

\noindent{\bf Proposition 3.1.} If $f\in L_1\cap L_2$, then $\|f\|_2=\| \Fo f\|_2=\|f \Fo \|_2$.

\

\noindent{\bf Proof}. Let $f'(x)=f^{\dag}(-x)$, one has
\begin{eqnarray*}
\|\Fo f\|_2^2
&=& \int_{\R^n}|\Fo f|^2dm(x) \\
&=& \int_{\R^n}<(\Fo f)^{\dag}(x)(\Fo f)(x)>_0dm(x) \\
&=& \int_{\R^n}<(f' \Fo)(x)(\Fo f)(x)>_0dm(x)\\
&=& <(f'*f)(0)>_0\\
&=& \int_{\R^n}<f'(-x)f(x)>_0dm(x) \\
&=& \int_{\R^n}<f^{\dag}(x)f(x)>_0dm(x)\\
&=& \|f\|_2^2 \\
\end{eqnarray*}
The proof for $f \Fo$ is analogous.

\hspace*{11cm} $\Box$

\

So the maps $f\to\Fo f$ and $f\to f\Fo$ can be uniquely extended in a continuous fashion to all of $L_2$. We will keep denoting this extension by $\Fo$. Moreover, these maps are surjective and satisfy Parseval's identity:

\

\noindent {\bf Parseval's Identity}. Let $f$ and $g\in L_2$, then
\begin{equation}\label{par1}
\int_{\R^n}f^{\dag}(x)g(x)dm(x)=\int_{\R^n}(\Fo f)^{\dag}(x)(\Fo g)(x)dm(x)
\end{equation}
and
\begin{equation}\label{par2}
\int_{\R^n}f(x)g^{\dag}(x)dm(x)=\int_{\R^n}(f\Fo )(x)(g\Fo )^{\dag}(x)dm(x).
\end{equation}

\

The scalar product in $L_2$ is defined using the scalar product of $\G_{2n}$ (\ref{prodesc}),
\begin{displaymath}
(f,g)_2=\int_{\R^n}f(x)\cdot g(x) \ dx.
\end{displaymath}
Then, equations (\ref{par1}) and (\ref{par2}) reduce to 
\begin{displaymath}
(f,g)_2=(\Fo f, \Fo g)_2=(f \Fo, g\Fo)_2
\end{displaymath}

\

\noindent {\bf Plancherel theorem.} The operators 
\begin{displaymath}
f\to \Fo f \quad \textrm{and} \quad f\to f\Fo
\end{displaymath}
from $L_2\to L_2$ are Hilbert space isomorphisms.

\

\section{Constructing Monogenic Functions}

\subsection{Monogenic Extensions using $\Fo$}

We will say that $f$ is monogenic (or left-monogenic) with respect to $B_i$ in a domain $M$ of $\R^{2n}$, if the (left) vector derivative of $f$ restricted to the subspace generated by $\{e_i,e'_i\}$ vanishes in all of $M$. That is,
\begin{displaymath}
\d_i f(x):= e_i\frac{\d f}{\d x_i}(x)+e'_i\frac{\d f}{\d x'_i}(x)=0, \quad \textrm{for all} \; x\in M.
\end{displaymath}
We define right-monogenic with respect to $B_i$ in $M$ in a similar way. If $f$ is monogenic (right-monogenic) with respect to each $B_i$, then $f$ is monogenic (right-monogenic), since 
\begin{displaymath}
\d=\d_1+\cdots+\d_n.
\end{displaymath}

The kernel of the Clifford-Fourier transform was defined using the symmetric bilinear function $I(x,y)=\sum_{i=1}^n B_iP_i(x)P_i(y)$. Since the projections $P_i$ are defined for all vectors in $\R^{2n}$, this map makes sense if we replace $y\in\R^n$ by $\yu\in\R^{2n}$, where
\begin{displaymath}
\yu=y+y'=(y_1e_1+\cdots+y_ne_n)+(y'_1e'_1+\cdots+y'_ne'_n)
\end{displaymath}
One has that
\begin{displaymath}
I(x,\yu)=I(x,y)-(x,y')
\end{displaymath}
where $(x,y')=x_1y'_1+\cdots +x_ny'_n$. 

Similarly it makes sense to talk about $I(\yu,x)$, but note that $I(x,\yu)\neq I(\yu,x)$ since
\begin{displaymath}
I(\yu,x)=I(x,y)+(x,y').
\end{displaymath}
Therefore, we can consider two extensions of the kernel:
\begin{displaymath}
e^{-I(x,\yu)} \qquad \textrm{and} \qquad e^{-I(\yu,x)}
\end{displaymath}

\noindent {\bf Lemma 4.1.} For each $x\in \R^n$, the functions $e^{\pm I(x,\yu)}$, $e^{\pm I(\yu,x)}\colon\R^{2n}\to\G_{2n}$, are monogenic and right-monogenic with respect to each $B_i$ in $\R^{2n}$, respectively.

\

Suppose $f\colon\R^n\to\G_{2n}$, assuming the corresponding integrals exist, we can define the following four extensions to $\R^{2n}$ of the Clifford-Fourier transform:

\begin{displaymath}
(\Fu f)(\yu):=\int_{\R^n}e^{-I(x,\yu)}f(x)dm(x),
\end{displaymath}
\begin{displaymath}
(f\Fu )(\yu):=\int_{\R^n}f(x)e^{-I(x,\yu)}dm(x),
\end{displaymath}
\begin{displaymath}
(\Fd f)(\yu):=\int_{\R^n}e^{-I(\yu,x)}f(x)dm(x),
\end{displaymath}
\begin{displaymath}
(f\Fd )(\yu):=\int_{\R^n}f(x)e^{-I(\yu,x)}dm(x).
\end{displaymath}

\

\noindent{\bf Proposition 4.1.} Let $f\colon\R^n\to\G_{2n}$ and $M$ be a domain of $\R^{2n}$.
\begin{enumerate}
\item If for each $\yu=y+y'\in M$ the functions $e^{(x,y')}f(x)$ and $x_i e^{(x,y')}f(x)$ are integrable, for $i=1,\dots,n$. Then, the extension of the Clifford-Fourier transform
\begin{displaymath}
(\Fu f)(\yu)
\end{displaymath}
is a monogenic function with respect to each $B_i$ in $M$.

\item If for each $\yu=y+y'\in M$ the functions $e^{-(x,y')}f(x)$ and $x_i e^{-(x,y')}f(x)$ are integrable, for $i=1,\dots,n$. Then, the extension of the Clifford-Fourier transform
\begin{displaymath}
(f \Fd)(\yu)
\end{displaymath}
is a right-monogenic function with respect to each $B_i$ en $M$.
\end{enumerate}

\

\noindent {\bf Proof.} Let $\yu\in M$, one has that
\begin{displaymath}
\frac{\d \Fu f}{\d y_i}(\yu)=-B_i (\Fu h)(\yu) \qquad \y \qquad \frac{\d \Fu f}{\d y'_i}(\yu)= (\Fu h)(\yu),
\end{displaymath}
where $h(x)= x_if(x)$. Since $e_iB_i=e'_i$, we have
\begin{displaymath}
\d_i (\Fu f)(\yu)=-e'_i (\Fu h)(\yu)+e'_i(\Fu h)(\yu)=0.
\end{displaymath}
The proof for $(f \Fd)(\yu)$ is analogous.

\hspace*{11cm} $\Box$

\

In many circumstances it is desirable to extend a continuous function in $\R^n$ to a monogenic or right-monogenic function in $\R^{2n}$. The following theorems present sufficient conditions for this to happen:

\

\noindent {\bf Theorem 4.1.} Let $f\colon\R^n\to\G_{2n}$ be continuous and integrable, suppose that the inverse Clifford-Fourier transform
\begin{displaymath}
(\Fo^{-1} f)(x)=\int_{\R^n}e^{I(x,t)}f(t)dm(t)
\end{displaymath} 
is integrable and for $\yu=y+y'\in\R^{2n}$, $e^{(x,y')}(\Fo^{-1} f)(x)$ and $x_ie^{(x,y')}(\Fo^{-1} f)(x)$ are integrable, for $i=1,\dots,n$. Then
\begin{displaymath}
(\Fu(\Fo^{-1} f))(\yu)
\end{displaymath}
is a monogenic extension of $f$ in all of $\R^{2n}$.

\

\noindent{\bf Proof.} Since $\Fo^{-1} f$ satisfies the conditions of proposition 4.1 with $M=\R^{2n}$, one has that $\Fu(\Fo^{-1}f)$ is monogenic in $\R^{2n}$. If $y\in \R^n$, then
\begin{displaymath}
(\Fu(\Fo^{-1}f))(y)=(\Fo(\Fo^{-1} f))(y)=f(y)
\end{displaymath}
since $f$ is continuous.

\hspace*{11cm} $\Box$

\

\noindent{\bf Theorem 4.2.} Let $f\colon\R^n\to\G_{2n}$ be continuous and integrable, assume that the inverse Clifford-Fourier transform
\begin{displaymath}
(f\Fo^{-1})(x)=\int_{\R^n}f(t)e^{I(x,t)}dm(t)
\end{displaymath} 
is integrable and that for all $\yu=y+y'\in\R^{2n}$, $e^{-(x,y')}(f\Fo^{-1})(x)$ and $ x_ie^{-(x,y')}(f\Fo^{-1})(x)$ are integrable, for $i=1,\dots,n$. Then the function
\begin{displaymath}
(( f\Fo^{-1})\Fd)(\yu)
\end{displaymath}
is a right-monogenic extension of $f$ in all of $\R^{2n}$.

\

Hence, the Clifford-Fourier transform gives useful machinery to extend a certain type of continuous functions to monogenic or right-monogenic functions. 

\

\subsection{Paley-Wiener Theorems}\label{teopw111}

In this section we will prove versions of the Paley-Wiener theorems for the Clifford-Fourier transform $\Fo$. These theorems show that under certain conditions the extension of the Clifford-Fourier transform extends a function in $L_2$ to a monogenic function. Conversely, they show that each monogenic function satisfying certain conditions, is the extension of the Clifford-Fourier transform of a function in $L_2$.

\

\noindent {\bf Theorem 4.3.}
\

\begin{enumerate}
\item Let $F\in L_2$ such that $F$ vanishes outside of $\R^n_+=(0,\infty)^n$ and let
\begin{equation}\label{teopw1}
f(\yu)=(\Fu F)(\yu), \qquad \yu\in\Pi^{2n}_-,
\end{equation}
where $\Pi^{2n}_-=\{\yu=y+y'\,|\, y \in \R^n, \, y'\in(-\infty,0)^n\}$. Then $f$ is monogenic with respect to each $B_i$, $i=1,\dots,n$, in $\Pi^{2n}_-$. Also, if $f_{y'}(y)=f(y+y')$ one has that
\begin{equation}\label{despw}
\sup_{y'\in (-\infty,0)^n} \int_{\R^n}|f_{y'}(y)|^2dm(y)= \|F\|_2^2.
\end{equation}
\item Conversely, if $f$ is monogenic with respect to each $B_i$, $i=1,\dots,n$, in $\Pi^{2n}_-$ and if there is a positive constant $C$ satisfying 
\begin{equation}\label{igupw}
\sup_{y'\in (-\infty,0)^n} \int_{\R^n}|f_{y'}(y)|^2dm(y)= C<\infty.
\end{equation}
Then there exists $F\in L_2$, vanishing outside of $\R^n_+$, such that (\ref{teopw1}) holds and $C=\|F\|_2^2$.
\end{enumerate}

\

\noindent{\bf Proof.} The first part of the theorem follows immediately using the results in section 3, Plancherel theorem, and the monotone convergence theorem.

For the second part we will use the notation 
\begin{displaymath}
\xk=x_1e_1+\cdots +x_{k-1}e_{k-1}+x_{k+1}e_{k+1}+\cdots+x_ne_n,
\end{displaymath}
that is, $\xk$ is the vector obtained after removing from $x$ the component $x_k$. Similarly we will use $\yk$, $\yyk$ and $\yuk=\yk+\yyk$. 

By lemma 4.1, one has that for each $x$ the function $e^{I(x,\yu)}f(\yu)$ is monogenic in $\Pi^{2n}_-$ with respect to each $B_i$. Now, fixing $k$, $y'_k<0$ ($y'_k\neq -1$) and $\alpha_k>0$, let 
\begin{displaymath}
E= \left\{
\begin{array}{ll}
[-1,y'_k], & y'_k > -1 \\
{[}y'_k,-1 {]},& y'_k < -1 \\
\end{array}
\right.
\end{displaymath}
By Clifford-Cauchy's theorem we get
\begin{eqnarray*}
0
&=& e_k \left( \int_{-\alpha_k}^{\alpha_k}e^{I_kx_k t}e^{-x_ky'_k}e^{I(\xk,\yuk)}f(te_k+y'_ke'_k+\yuk)dm(t)\right. \\
&-&  \left. \int_{-\alpha_k}^{\alpha_k}e^{I_kx_k t}e^{-x_k(-1)}e^{I(\xk,\yuk)}f(te_k+(-1)e'_k+\yuk)dm(t) \right) \\
&+& e'_k \left(\int_{E}e^{I_kx_k\alpha_k}e^{-x_kt}e^{I(\xk,\yuk)}f(\alpha_ke_k+te'_k+\yuk)dm(t)\right.\\
&-& \left. \int_{E}e^{I_kx_k(-\alpha_k)}e^{-x_kt}e^{I(\xk,\yuk)}f((-\alpha_k)e_k+te'_k+\yuk)dm(t)\right)
\end{eqnarray*}
By Fubini's theorem and condition (\ref{igupw}), we get a sequence $\{\alpha_{k,j}\}_{j}$ such that $\alpha_{k,j}\to\infty$ and
\begin{eqnarray*}
0
&=& \lim_{j\to\infty} \left[ \int_{-\alpha_{k,j}}^{\alpha_{k,j}}e^{I_kx_k t}e^{-x_ky'_k}e^{I(\xk,\yuk)}f(te_k+y'_ke'_k+\yuk)dm(t) \right.\\
&-& \left. \int_{-\alpha_{k,j}}^{\alpha_{k,j}}e^{I_kx_k t}e^{-x_k(-1)}e^{I(\xk,\yuk)}f(te_k+(-1)e'_k+\yuk)dm(t)\right]\\
\end{eqnarray*}
for almost all $\yuk$. We can do this for each $k=1,\dots,n$. Hence, if $A_j=[-\alpha_{1,j},\alpha_{1,j}]\times\cdots\times[-\alpha_{n,j},\alpha_{n,j}]$ and $\u=e'_1+\cdots +e'_n$, then
\begin{equation}\label{equa}
\lim_{j\to\infty}\left[ \int_{A_j}e^{I(x,y+y')}f(y+y')dm(y)-\int_{A_j}e^{I(x,y-\u)}f(y-\u)dm(y) \right]=0
\end{equation}
for all $x\in\R^n$. Using the notation $f_{-\u}(y)=f(y-\u)$, define 
\begin{displaymath}
F(x)=e^{-(x,-\u)}(\Fo^{-1}f_{-\u})(x)
\end{displaymath}
by Plancherel's theorem and (\ref{equa}) we get
\begin{displaymath}
F(x)=e^{-(x,y')}(\Fo^{-1}f_{y'})(x), \; \mathrm{almost \ everywhere, \ for }\ y'\in (-\infty,0)^n.
\end{displaymath}
Again, by Plancherel's theorem, we get
\begin{equation}\label{CF2}
\int_{\R^n}e^{2(x,y')}|F(x)|^2dm(x)=\int_{\R^n}|f_{y'}(x)|^2dm(x)\leq C.
\end{equation}
This shows that for $x\notin\R^n_+$, if $y'=-\lambda\u$ and $\lambda\to \infty$, then $F(x)=0$ almost everywhere. Using the monotone convergence theorem, this shows that if $y'\to 0$ in $(-\infty,0)^n$ one has that $\|F\|^2_2\leq C$. Therefore, $F$ is a function in $L_2$ vanishing outside $\R^n_+$.

Finally, since $F$ vanishes outside $\R^n_+$ and $(\Fo^{-1}f_{y'})(x)=e^{(x,y')}F(x)$ almost everywhere, then $(\Fo^{-1}f_{y'})\in L_1$ for $y'\in(-\infty,0)^n$. Hence, for $\yu\in\Pi^{2n}_-$ one has that
\begin{displaymath}
f(\yu)=f_{y'}(y)=\int_{\R^n}e^{-I(x,y)}(\Fo^{-1}f_{y'})(x)dm(x)=(\Fu F)(\yu),
\end{displaymath}
since both functions are continuous. By the last equation we have $\|f_{y'}\|_2^2\leq\|F\|_2^2$ and in (\ref{CF2}) we found that $\|F\|_2^2\leq C$, putting these inequalities together we get $C=\|F\|_2^2$.

\hspace*{11cm} $\Box$

\

\noindent {\bf Remark.} One can show an analogous result for right-monogenic functions using the extension $(F \, \Fd)(\yu)$ of the Clifford-Fourier transform.

\

\noindent {\bf Theorem 4.4.}
\

\begin{enumerate}
\item Let $A$ be a positive constant and $B_A$ the ball in $\R^n$ centered at zero of radius $A$. Let $F\in L_2$ vanishing outside of $B_A$ and define
\begin{equation}\label{teopw21}
f(\yu)=(\Fu F)(\yu),\qquad \yu\in\R^{2n}.
\end{equation}
Then $f$ is monogenic with respect to each $B_i$, $i=1,\dots,n$, in all of $\R^{2n}$. Also, there is a positive constant $C$ such that
\begin{equation}\label{teopw22}
|f(\yu)|\leq C\,e^{A|\yu|}
\end{equation}
and if $f_{y'}(y)=f(y+y')$ one has
\begin{equation}\label{teopw23}
\|f_{y'}\|_2\leq e^{A|y'|}\|F\|_2.
\end{equation}
\item Conversely, if $A$ and $C$ are positive constants and $f$ is monogenic with respect to each $B_i$, $i=1,\dots,n$, in $\R^{2n}$ satisfying inequality (\ref{teopw22}) and
\begin{equation}\label{cond21}
\|f_{y'}\|_2^2\leq h(y')<\infty,
\end{equation}
where $h$ is a locally integrable function.  Then there is $F\in L_2$ vanishing outside $B_A$ such that (\ref{teopw21}) holds.
\end{enumerate}

\

\noindent {\bf Proof.} The first part follows  by the results in section 3 and Plancherel's theorem.

For the second part we use Fubini's theorem, condition (\ref{cond21}) and Clifford-Cauchy's theorem to get, as in theorem 4.3,  a sequence $\{\alpha_{k,j}\}_j$ such that $\alpha_{k,j}\to\infty$ and
\begin{eqnarray*}
0
&=& \lim_{j\to\infty} \left[ \int_{-\alpha_{k,j}}^{\alpha_{k,j}}e^{I_kx_k t}e^{-x_ky'_k}e^{I(\xk,\yuk)}f(te_k+y'_ke'_k+\yuk)dm(t) \right. \\
&-& \left. \int_{-\alpha_{k,j}}^{\alpha_{k,j}}e^{I_kx_k t}e^{-x_k \yb'_k}e^{I(\xk,\yuk)}f(te_k+\yb'_ke'_k+\yuk)dm(t)\right],
\end{eqnarray*}
for almost all $\yuk$ and $y'_k\neq\yb'_k$. Applying this for each $k=1.\dots,n$, we get 
\begin{displaymath}
\lim_{j\to\infty}\left[ \int_{A_j}e^{I(x,y+y')}f(y+y')dm(y)-\int_{A_j}e^{I(x,y+\yb')}f(y+\yb')dm(y) \right]=0,
\end{displaymath}
for each $x\in\R^n$, where $A_j=[-\alpha_{1,j},\alpha_{1,j}]\times\cdots\times[-\alpha_{n,j},\alpha_{n,j}]$. The previous equation is valid for almost all $y'$ and $\yb'$. We fix $\yb'$ such that this expression holds. Define
\begin{displaymath}
F(x)=e^{-(x,\yb')}(\Fo^{-1}f_{\yb'})(x)
\end{displaymath}
by Plancherel's theorem, we get 
\begin{displaymath}
F(x)=e^{-(x,y')}(\Fo^{-1}f_{y'})(x),\quad \mathrm{for \ almost \ all} \ x \ \y \ y'.
\end{displaymath}
Assume $x\neq 0$, for $\lambda>0$ and any $\epsilon>0$ we can find $u=u(\lambda,\epsilon)\in \R^n$ such that $|u|<\epsilon$ and 
\begin{displaymath}
y'=(\frac{\lambda x_1}{|x|}+u_1)e'_1+\cdots+(\frac{\lambda x_n}{|x|}+u_n)e'_n\;\in\;M,
\end{displaymath} 
then $|y'|<\lambda+\epsilon \ $ and $ \ (x,y')=\lambda|x|+(x,u)$. Thus, for each $r>0$, 
\begin{displaymath}
\left|e^{-(x,y')}\int_{B_r}e^{I(x,y)}f_{y'}(y)dm(y)\right|\leq 2^n\, e^{-\lambda|x|}e^{-(x,u)}\int_{B_r}Ce^{A|\yu|}dm(y)
\end{displaymath}
\begin{displaymath}
\leq \left(2^n\,C \int_{B_r}e^{A|y|}dm(y)\right) e^{(A\epsilon-(x,u))}e^{\lambda(A-|x|)}.
\end{displaymath}
Now, if $\epsilon\to 0$ then $e^{(A\epsilon-(x,u))}\to1$ and if $\lambda\to\infty$, when $|x|>A$, we get $e^{\lambda(A-|x|)}\to 0$. 

By Plancherel's theorem we get $e^{-(x,y')}(\Fo^{-1}f_{y'})(x)=0$ almost everywhere outside of $B_A$. Hence, $F(x)=0$ almost everywhere outside of $B_A$, and then, by definition, $F\in L_2$. Thus, since $(\Fo^{-1} f_{y'})(x)=e^{(x,y')}F(x)$ for almost all $x$ and $y'$, we get $(\Fo^{-1} f_{y'})\in L_1$. Therefore, for $\yu\in\R^{2n}$ we get
\begin{displaymath}
f(\yu)=f_{y'}(y)=\int_{\R^n}e^{-I(x,y)}(\Fo^{-1}f_{y'})(x)dm(x)=(\Fu F)(\yu),
\end{displaymath}
since both functions are continuous.

\hspace*{11cm} $\Box$

\

\noindent {\bf Remark.} Similarly, one can show an analogous result for right-monogenic functions using the extension $(F \, \Fd)(\yu)$.

\label{end-art}

\end{document}